\documentclass[reqno,12pt]{amsart}
\theoremstyle{plain}    
\newtheorem{thm}{Theorem}[section]
\theoremstyle{plain}    
\newtheorem{lemma}[thm]{Lemma} 
\newtheorem{prop}[thm]{Proposition}
\theoremstyle{remark}
\newtheorem{defi}[thm]{Definition}

\newtheorem{remark}[thm]{Remark}

\usepackage{amscd,amssymb,comment,epic,euscript,graphics}

\newcommand\Afr{{\mathfrak A}}

\newcommand\alphat{{\tilde\alpha}}

\newcommand\Bc{{\mathcal{B}}}

\newcommand\betat{{\tilde\beta}}

\newcommand\Cc{{\mathcal{C}}}

\newcommand\cinv{^{\langle-1\rangle}}

\newcommand\Cpx{{\mathbf C}}

\newcommand\Ec{{\mathcal{E}}}

\newcommand\Fc{{\mathcal{F}}}

\newcommand\id{{\operatorname{id}}}

\newcommand\Nats{{\mathbf N}}

\newcommand\op{{\rm op}}

\newcommand\otdt{\otimes\cdots\otimes}

\newcommand\oth{\hat{\otimes}}

\newcommand\tdt{\times\cdots\times}

\begin{document}

 \title{On the S--transform over a Banach algebra}

 \author{Kenneth J.\ Dykema}

 \address{\hskip-\parindent
 Mathematisches Institut \\
 Westf\"alische Wilhelms--Universit\"at M\"unster \\
 Einsteinstr.\ 62 \\
 48149 M\"unster \\
 Germany}

 \address{\hskip-\parindent
 {\bf Permanent Address:} Department of Mathematics \\
 Texas A\&M University \\
 College Station TX 77843--3368, USA}
 \email{kdykema@math.tamu.edu}

 \thanks{Supported in part by NSF grant DMS--0300336 and by the Alexander von Humboldt Foundation.}

 \date{11 January, 2005}

\begin{abstract}
The S--transform is shown to satisfy a specific
twisted multiplicativity property
for free random variables
in a $B$--valued Banach noncommutative probability space, for an arbitrary
unital complex Banach algebra $B$.
Also, a new proof of the additivity of the R--transform in this setting is given.
\end{abstract}

 \maketitle


\section{Introduction and statement of the main result}

Let $B$ be a unital complex Banach algebra.
(In this paper, all Banach algebras will be over the complex numbers.)
A $B$--valued Banach noncommutative probability space
is a pair $(A,E)$ where $A$ is a unital Banach algebra containing
an isometrically embedded copy of $B$
as a unital subalgebra and where $E:A\to B$ is a bounded projection
satisfying the conditional expectation property
\[
E(b_1ab_2)=b_1E(a)b_2\qquad(a\in A,\,b_1,b_2\in B).
\]
In the free probability theory of Voiculescu, see~\cite{V85} and~\cite{VDN},
elements $x$ and $y$ of $A$ are said to be free if their mixed moments
$E(b_1a_1\cdots b_na_n)$, where $a_j\in\{x,y\}$ and $b_j\in B$, are determined
in a specific way from the moments of $x$ and of $y$.
Of particular interest, for example to garner spectral data, are the symmetric
moments 
\begin{equation}\label{eq:bxy}
E(bxybxy\cdots bxy)
\end{equation}
of the product $xy$, for $b\in B$.

In the case $B=\Cpx$, Voiculescu~\cite{V87} invented the S--transform of
an element $x\in A$ satisfying $E(x)\ne0$.
The S--transform can be used to find the generating function for
the symmetric moments~\eqref{eq:bxy}
of $xy$ in terms of those for $x$ and $y$ individually, when $x$ and $y$ are free
and when $E(x)\ne0$ and $E(y)\ne0$.
In particular, Voiculescu showed that the S--transform
is multiplicative:
\begin{equation}\label{eq:Smult}
S_{xy}=S_xS_y
\end{equation}
when $x$ and $y$ are free.

In~\cite{V95}, Voiculescu gave a definition of an S--transform in the context of an arbitrary
noncommutative probability space.
However, this definition was quite complicated and involved differential equations.

Recently, Aagaard~\cite{Aa} took the straightforward extension
of Voiculescu's definition~\cite{V87} of the scalar--valued S--transform
to the Banach algebra situation
and generalized Voiculescu's result~\eqref{eq:Smult}
to the case when $B$ is a commutative unital Banach algebra
and $E(x)$ and $E(y)$ are invertible elements of $B$.

In this paper, we treat the case when $B$ is an arbitrary unital Banach algebra.
We make an improvement in Aagaard's definition
of the S--transform.
For us, $S_x$ is a $B$--valued analytic function defined in a neighborhood of $0$ in $B$.
We write $S_{xy}$ in terms of $S_x$ and $S_y$
(again assuming $E(x)$ and $E(y)$ are invertible).
Instead of simple multiplicativity~\eqref{eq:Smult}, we have in general
a twisted multiplicativity, as stated in our main theorem immediately below,
which reduces to~\eqref{eq:Smult} when $B$ is commutative.

\begin{thm}\label{thm:main}
Let $B$ be a unital complex Banach algebra and let $(A,E)$ be a $B$--valued
Banach noncommutative probability space.
Let $x,y\in A$ be free in $(A,E)$ and assume both $E(x)$ and $E(y)$ are invertible
elements of $B$.
Then
\begin{equation}\label{eq:main}
S_{xy}(b)=S_y(b)S_x(S_y(b)^{-1}bS_y(b)).
\end{equation}
\end{thm}

Our definition of the S--transform and our proof of Theorem~\ref{thm:main}
rely on the theory of analytic functions between Banach spaces --
see for example Chapters~III and~XXVI of~\cite{HP} and papers cited there.

In~\cite{H}, Haagerup gave two new proofs of the multiplicativity of the S--transform
in the case $B=\Cpx$.
Our proof of Theorem~\ref{thm:main} is very much inspired by one of Haagerup's proofs,
namely Theorem~2.3 of~\cite{H}, which 
uses creation and annihilation operators in the full Fock space.
In particular, we consider a $B$--valued Banach algebra
analogue of the full Fock space and we construct random variables having arbitrary moments
up to a given finite order, using analogues of the creation and annihilation
operators.
These are reminiscent of, though slightly different from, Voiculescu's constructions in~\cite{V95}.

In~\S\ref{sec:def} below, we define the S--transform $S_a$
(assuming the expectation of $a$ is invertible).
Then, considering Taylor expansions about zero, we show that the $n$th order term in the expansion
for $S_a$ depends only on the moments up to $n$th order of $a$.
In~\S\ref{sec:tmult}, we construct operators
analogous to the creation and annihilation operators on full Fock space, and we use these to
prove the main result, Theorem~\ref{thm:main}.
In~\S\ref{sec:R-t}, we offer a new proof of additivity of the R--tranfrom over
a Banach space, using the operators and techniques introduced in the preceding sections.

\smallskip
\noindent{\em Acknowledgements.} The author wishes to thank
Joachim Cuntz and the Mathematics Institute of the
Westf\"alische Wilhelms--Universit\"at
M\"unster for their generous hospitality during the author's year--long visit, when
this research was conducted.

\section{The S--transform in a Banach noncommutative probability space}
\label{sec:def}

Let $B$ be a unital Banach algebra.
For $n\ge1$ we will let $\Bc_n(B)$ denote the set of all bounded $n$--multilinear maps
\[
\alpha_n:\underset{n\text{ times}}{\underbrace{B\tdt B}}\to B,
\]
where multilinearity means over $\Cpx$ and a multilinear map $\alpha_n$ is bounded if
\[
\|\alpha_n\|:=\sup\{\|\alpha_n(b_1,\ldots,b_n)\|\mid b_j\in B,\,\|b_1\|,\ldots,\|b_n\|\le1\}<\infty.
\]
We say $\alpha_n$ is {\em symmetric} if it is invariant under arbitrary
permutations of its $n$ arguments.

From the theory of analytic functions between complex Banach spaces,
any $B$--valued
analytic function $F$ defined on a neighborhood of zero in $B$ has an
expansion
\begin{equation}\label{eq:Fexp}
F(b)=F(0)+\sum_{n=1}^\infty F_n(b,\ldots,b),
\end{equation}
for some symmetric multilinear functions $F_n\in\Bc_n(B)$,
with $\limsup_{n\to\infty}\|F_n\|^{1/n}<\infty$;
see, for example, Theorem 3.17.1 of~\cite{HP} and its proof.
Here $F_1$ is just the Fr\'echet derivative of $F$ at $0$
and the multilinear function $F_n$ appearing in~\eqref{eq:Fexp}
is $1/n!$ times the $n$th variation of $F$,
i.e.\ $n!F_n(h_1,\ldots,h_n)$ is the $n$--fold Fr\'echet derivative
taken with respect to increments $h_1,\ldots,h_n$.
For convenience we will write $F_0$ for $F(0)$.
We will refer to~\eqref{eq:Fexp} as the
{\em power series expansion} of $F(b)$ around $0$ and to $F_n(b,\ldots,b)$ as the
$n$th term in this power series expansion.
Note that the full symmetric multilinear function $F_n$ can be recovered from
knowing its diagonal $b\mapsto F_n(b,\ldots,b)$;
for example, $n!F_n(b_1,\ldots,b_n)$ is the obvious partial derivative of
\[
F_n(t_1b_1+\cdots+t_nb_n,\ldots,t_1b_1+\cdots+t_nb_n)
\]
at $(0,\ldots,0)$, where $t_1,\ldots,t_n$
are real variables.

Let $(A,E)$ be a Banach
noncommutative probability space over $B$, let $a\in A$ and suppose $E(a)$
is an invertible element of $B$.
Consider the function
\begin{equation}\label{eq:Psi}
\Psi_a(b)=E((1-ba)^{-1})-1=\sum_{n=1}^\infty E((ba)^n),
\end{equation}
defined for $\|b\|<\|a\|^{-1}$.
Then $\Psi_a$ is Fr\'echet differentiable on its domain, i.e.\ is analytic there.
We also have 
\begin{equation}\label{eq:PsiPhi}
\Psi_a(b)=b\Phi_a(b),
\end{equation}
where 
\begin{equation}\label{eq:Phi}
\Phi_a(b)=E(a(1-ba)^{-1});
\end{equation}
clearly $\Phi_a$
is analytic on the domain of $\Psi_a$.
The Fr\'echet differential of $\Psi_a$ at $b=0$ is easily found to be the bounded linear map
\begin{equation}\label{eq:hE}
h\mapsto hE(a)
\end{equation}
from $B$ to itself.
By hypothesis, this linear map has bounded inverse $h\mapsto hE(a)^{-1}$.
By the usual Banach space inverse function theorem, there are neighborhoods $U$ and $V$
of zero in $B$ such that $U$ lies in the domain of $\Psi_a$ and the restriction of $\Psi_a$
to $U$ is a homeomorphism onto $V$.
Moreover, letting $\Psi_a\cinv$ denote the inverse with respect to composition
of the restriction of $\Psi_a$ to $U$,
the function $\Psi_a\cinv$ is Fr\'echet differentiable on its domain
and is, therefore, analytic there.

\begin{lemma}\label{lem:Ha}
Assuming $E(a)$ is invertible,
there is an open neighborhood of $0$ in $B$ and
unique analytic $B$--valued function $H_a$
defined there
such that $\Psi_a\cinv(b)=bH_a(b)$.
\end{lemma}
\begin{proof}
Uniqueness of $H_a$ is clear by uniqueness of power series expansions about zero.
Let us show existence.
Using~\eqref{eq:PsiPhi}, we seek $H_a$ such that $bH_a(b)\Phi_a(bH_a(b))=b$, and it will suffice to
find $H_a$ such that
\begin{equation}\label{eq:HPhi}
H_a(b)\Phi_a(bH_a(b))=1.
\end{equation}
The existence of $H_a$ follows from
an easy application of the implicit function theorem for functions 
between Banach spaces, which is a result of Hildebrandt and Graves~\cite{HG}
(see also the discussion on p.\ 655 of~\cite{G}).
Indeed, $H_a(0)=E(a)^{-1}$ is a solution of~\eqref{eq:HPhi} at $b=0$
and the Fr\'echet differential of the function $x\mapsto x\Phi_a(bx)$
at $b=0$ is the map~\eqref{eq:hE},
which has bounded inverse.
\end{proof}

\begin{defi}
Let $a\in A$ and assume $E(a)$ is invertible.
The {\em S--transform} of $a$ is the $B$--valued analytic function
\begin{equation}\label{eq:S}
S_a(b)=(1+b)H_a(b),
\end{equation}
which defined in some neighborhood of $0$ in $B$,
where $H_a$ is the function from Lemma~\ref{lem:Ha}.
\end{defi}
Note that $S_a(0)=E(a)^{-1}$.

We may write 
\begin{equation}\label{eq:AaS}
S_a(b)=(1+b)b^{-1}\Psi_a\cinv(b),
\end{equation}
which is the same formula given by Voiculescu~\cite{V87} and used by Aagaard~\cite{Aa}.
In the case $B=\Cpx$, the definition~\eqref{eq:S} yields, of course,
the same function as Voiculescu's S--transform.
Moreover, the only difference between the definition~\eqref{eq:S} and the one appearing in~\cite{Aa}
is that we have used the implicit function theorem to show that~\eqref{eq:AaS}
makes sense for all $b$ in a neighborhood of zero.

If $F$, $G$ and $H$ are $B$--valued analytic functions defined on neighborhoods of $0$ in $B$,
then the product $FG$ is analytic
and, if $H(0)=0$, also the composition $F\circ H$ is analytic in some neighborhood of $0$ in $B$.
Straightforward asymptotic analysis yields the following formulas for the
diagonals of the multilinear functions
appearing in the power series expansions of $FG$ and $F\circ H$.
\begin{lemma}\label{lem:FG}
We have for $n\ge0$
\begin{equation}\label{eq:FG}
(FG)_n(b,\ldots,b)=\sum_{k=0}^nF_k(b,\ldots,b)G_{n-k}(b,\ldots,b)
\end{equation}
and for $n\ge1$
\begin{equation}\label{eq:FoH}
(F\circ H)_n(b,\ldots,b)=
\sum_{k=1}^n\sum_{\substack{p_1,\ldots,p_k\ge1 \\ p_1+\cdots+p_k=n}}
F_k(H_{p_1}(b,\ldots,b),\ldots,H_{p_k}(b,\ldots,b)).
\end{equation}
\end{lemma}
\begin{lemma}\label{lem:Finv}
Let $F$ be analytic in a neighborhood of $0$.
If $F(0)$ is an invertible element of $B$, then $G(b)=F(b)^{-1}$ defines
a function that is analytic in a neighborhood of $0$, and the $n$th term of
its power series expansion is $G_0=F_0^{-1}$ and, for $n\ge1$,
\begin{equation}\label{eq:Frecip}
G_n(b,\ldots,b)=-F_0^{-1}\sum_{k=1}^nF_k(b,\ldots,b)G_{n-k}(b,\ldots,b).
\end{equation}
On the other hand, if $F(0)=0$ and if $F_1$ has a bounded inverse, then
$F$ has an inverse with respect to composition, denoted $F\cinv$, that is analytic
in a neighborhood of $0$.
Taking $H=F\cinv$, we have
$H_1=(F_1)\cinv$ and, for $n\ge2$,
\begin{equation}\label{eq:Finv}
H_n(b,\ldots,b)=
-(F_1)\cinv\bigg(\sum_{k=2}^n\sum_{\substack{p_1,\ldots,p_k\ge1 \\ p_1+\cdots+p_k=n}}
F_k(H_{p_1}(b,\ldots,b),\ldots,
H_{p_k}(b,\ldots,b))\bigg).
\end{equation}
\end{lemma}
\begin{proof}
Assuming $F(0)$ is invertible, that $G(b)=F(b)^{-1}$ is analytic is clear,
and we have $(FG)_0=1$ and $(FG)_n=0$ for $n\ge1$.
Now the expression~\eqref{eq:Frecip} results from solving~\eqref{eq:FG} for $G_n$.

If $F(0)=0$ and the Fr\'echet derivative $F_1$ of $F$ at $0$ has bounded inverse, then by the inverse
function theorem for Banach spaces,
$F$ has an inverse with respect to composition $F\cinv$ that is analytic
in a neighborhood of $0$.
Taking $H=F\cinv$, we have $(F\circ H)_1=\id_B$ and $(F\circ H)_n=0$ for all $n\ge2$.
Solving in~\eqref{eq:FoH} for $H_n$ yields the expression~\eqref{eq:Finv}.
\end{proof}

Consider an element $a\in A$ as
at the begining of this section.
We say the {\em $n$th moment function} of $a$ is the multilinear function $\mu_{a,n}\in\Bc_n(B)$
given by
\[
\mu_{a,n}(b_1,\ldots,b_n)=E(b_1ab_2a\cdots b_na).
\]
\begin{prop}\label{prop:nth}
Assume $E(a)$ is an invertible element of $B$.
Then the
$n$th term $(S_{a})_n(b,\ldots,b)$
in the power series expansion of the S--transform $S_a$ of $a$
about zero
depends only on the first $n$ moment functions
$\mu_{a,1},\,\mu_{a,2},\,\ldots,\,\mu_{a,n}$ of $a$.
\end{prop}
\begin{proof}
The symmetric $n$--multilinear function $(\Psi_a)_n$
appearing in the power series expansion of $\Psi_a$ is the symmetrization of $\mu_{a,n}$.
Using Lemma~\ref{lem:Finv}, we see that the $n$th term $(\Psi_a\cinv)_n(b,\ldots,b)$
in the power series expansion of $\Psi_a\cinv(b)$ around $0$ depends only on
$\mu_{a,1},\ldots,\mu_{a,n}$.
But 
\begin{gather*}
(\Psi_a\cinv)_n(b,\ldots,b)=b\,(H_a)_n(b,\ldots,b) \\
(S_{a})_n(b,\ldots,b)=(1+b)\,(H_a)_n(b,\ldots,b)
\end{gather*}
and the result is proved.
\end{proof}

\section{Twisted multiplicativity of the S--transform}
\label{sec:tmult}

Let $B$ be a unital Banach algebra over $\Cpx$ and let $I$ be a set.
Let $D=\ell^1(I,B)$ be the Banach space of all functions $d:I\to B$
such that $\|d\|:=\sum_{i\in I}\|d(i)\|<\infty$.
For $i\in I$, $\delta_i\in D$ will denote the function taking value $1$ at $i$
and $0$ at all other elements of $I$.
We have the obvious left action of $B$ on $D$ by $(bd)(i)=b\,d(i)$,
and the resulting algebra homomorphism $B\to\Bc(D)$ is isometric.
(Whenever $X$ is a Banach space, we denote by $\Bc(X)$ the Banach algebra
of all bounded linear operators from $X$ to itself.)
For $k\ge1$, let $D^{\oth k}=D\oth\cdots\oth D$ be the $k$--fold Banach space projective
tensor product of $D$ with itself (over the complex field).
Consider the Banach space
\begin{equation}\label{eq:F}
\Fc=B\Omega\oplus\bigoplus_{k=1}^\infty D^{\oth k}\oth B,
\end{equation}
where also $\oth B$ is the Banach space projective tensor product and where the we take
the direct sum with respect to the $\ell^1$--norm.
Here, $B\Omega$ signifies just a copy of $B$ and $\Omega$ denotes the identity element of this
copy of $B$, consdered as a vector in $\Fc$.
Let $\lambda:B\to\Bc(\Fc)$ be the map defined by
\begin{gather*}
\lambda(b)(b_0\Omega)=(bb_0)\Omega \\
\lambda(b)(d_1\otdt d_k\otimes b_0)=(bd_1)\otimes d_2\otdt d_k\otimes b_0
\end{gather*}
for $k\in\Nats$, $d_1,\ldots,d_k\in D$ and $b_0\in B$.
Then $\lambda$ is an isometric algebra homomorphism.
We will often omit to write $\lambda$, and just think of $B$ as included in $\Bc(\Fc)$
by this left action.

\begin{remark}
For specificity, we took the $\ell^1$ norms in the definitions of $D$ and $\Fc$,
but we actually have considerable flexibility.
For $D$ we need only a Banach space completion of the set of all functions
$d:I\to B$ vanishing at all but finitely many elements in $I$ with the property
$\|b\delta_i\|=\|b\|$, and similarly for $\Fc$.
Moreover, we could replace the projective tensor norm $\oth B$ in~\eqref{eq:F}
with any tensor norm so that $\|d\otimes B\|=\|d\|\,\|b\|$ for all $d\in D^{\oth k}$
and $b\in B$.
\end{remark}

Let $P:\Fc\to B$ be the projection onto the summand $B\Omega=B$
that sends all summands $D^{\oth k}\oth B$ to zero
and let $\Ec:\Bc(\Fc)\to B$
be $\Ec(X)=P(X\Omega)$.
Then $\Ec$ has norm $1$ and satisfies $\Ec\circ\lambda=\id_B$.
Let $\rho:B\to\Bc(\Fc)$ be the map defined by
\begin{gather*}
\rho(b)(b_0\Omega)=(b_0b)\Omega \\
\rho(b)(d_1\otdt d_k\otimes b_0)= d_1\otdt d_k\otimes (b_0b).
\end{gather*}
Then $\rho$ is an isometric algebra isomorphism from the opposite algebra
$B^\op$ into $\Bc(\Fc)$.
Let $\Bc(\Fc)\cap\rho(B)'$ denote the set of all bounded operators on $\Fc$
that commute with $\rho(b)$ for all $b\in B$.
Note that $\lambda(B)\subseteq\Bc(\Fc)\cap\rho(B)'$.
\begin{prop}
The restriction of $\Ec$ to $\Bc(\Fc)\cap\rho(B)'$ satisfies the conditional expectation
property
\[
\Ec(b_1Xb_2)=b_1\Ec(X)b_2\qquad(X\in\Bc(\Fc)\cap\rho(B)',\,b_1,b_2\in B).
\]
\end{prop}
\begin{proof}
We have
\begin{align*}
\Ec(b_1Xb_2)
&=P(\lambda(b_1)X\lambda(b_2)\Omega)
=P(\lambda(b_1)X\rho(b_2)\Omega) \\
&=P(\rho(b_2)\lambda(b_1)X\Omega)
=P(\lambda(b_1)X\Omega)b_2
=b_1P(X\Omega)b_2
=b_1\Ec(X)b_2.
\end{align*}
\end{proof}

For $i\in I$, let $L_i\in\Bc(\Fc)$ be defined by
\begin{align*}
L_i(b_0\Omega)&=\delta_i\otimes b_0 \\
L_i(d_1\otdt d_k\otimes b_0)&=\delta_i\otimes d_1\otdt d_k\otimes b_0.
\end{align*} 
Thus,
\[
b_1\delta_{i_1}\otimes b_2\delta_{i_2}\otdt b_k\delta_{i_k}\otimes b_0
=b_1L_{i_1}b_2L_{i_2}\cdots b_kL_{i_k}b_0\Omega.
\]
Recall that $\Bc_n(B)$ denotes the set of all bounded multilinear functions from the $n$--fold
product of $B$ to $B$.
We will also let $\Bc_0(B)=B$.
If $i\in I$, $n\in\Nats$ and
$\alpha_n\in\Bc_n(B)$, define $V_{i,n}(\alpha_n)$ and $W_{i,n}(\alpha_n)$ in $\Bc(\Fc)$
by
\begin{align*}
V_{i,n}(\alpha_n)(b_0\Omega)&=0 \\[1ex]
V_{i,n}(\alpha_n)(d_1\otdt d_k\otimes b_0)&=
\begin{cases}
0,&k<n \\
\alpha_n(d_1(i),\ldots,d_n(i))b_0\Omega,&k=n \\
\alpha_n(d_1(i),\ldots,d_n(i))d_{n+1}\otdt d_k\otimes b_0,&k>n
\end{cases}
\end{align*}
and
\begin{align*}
W_{i,n}(\alpha_n)&(b_0\Omega)=0 \displaybreak[2] \\[1ex]
W_{i,n}(\alpha_n)&(d_1\otdt d_k\otimes b_0)= \\
&=\begin{cases}
0,&k<n \\
\alpha_n(d_1(i),\ldots,d_n(i))\delta_i\otimes b_0,&k=n \\
\alpha_n(d_1(i),\ldots,d_n(i))\delta_i\otimes d_{n+1}\otdt d_k\otimes b_0,&k>n.
\end{cases}
\end{align*}
Finally, taking $n=0$ and $\alpha_0\in B$, let
\[
V_{i,0}(\alpha_0)=\alpha_0\qquad W_{i,0}(\alpha_0)=\alpha_0L_i.
\]
These formulas are guaranteed to define bounded operators on $\Fc$, because we
took the projective tensor product in $D^{\oth k}$.
The expression $V_{i,n}(\alpha_n)$, $n\ge1$, is a sort of $n$--fold annihilation operator,
while $W_{i,n}(\alpha_n)$ is $n$--fold annihilation combined with single creation, and, of course,
$W_{i,0}$ is a single creation operator.
Note that in all cases we have $V_{i,n}(\alpha_n),\,W_{i,n}(\alpha_n)\in\Bc(\Fc)\cap\rho(B)'$.

The relations gathered in the following lemma are easily verified.
\begin{lemma}\label{lem:rel}
Let $n,m\in\Nats$ and
$\alpha_n\in\Bc_n(B)$,
$\beta_m\in\Bc_{m}(B)$
and take $b\in B$.
Then
\renewcommand{\labelenumi}{(\roman{enumi})}
\begin{enumerate}
\item
\[
V_{i,n}(\alpha_n)\lambda(b)=V_{i,n}(\alphat_n)\qquad W_{i,n}(\alpha_n)\lambda(b)=W_{i,n}(\alphat_n),
\]
where 
\[
\alphat_n(b_1,\ldots,b_n)=\alpha_n(bb_1,b_2,\ldots,b_n);
\]
\item if $n=1$, then
\[
V_{i,1}(\alpha_1)L_i=\lambda(\alpha_1(1)),\qquad W_{i,1}(\alpha_1)L_i=\lambda(\alpha_1(1))L_i
\]
and for $n\ge2$ we have
\[
V_{i,n}(\alpha_n)L_i=V_{i,n-1}(\alphat_{n-1}),\qquad W_{i,n}(\alpha_n)L_i=W_{i,n-1}(\alphat_{n-1}),
\]
where here 
\[
\alphat_{n-1}(b_1,\ldots,b_{n-1})=\alpha_n(1,b_1,\ldots,b_{n-1});
\]
\item we have
\[
V_{i,n}(\alpha_n)V_{i,m}(\beta_m)=V_{i,n+m}(\gamma_{n+m}),
\qquad W_{i,n}(\alpha_n)V_{i,m}(\beta_m)=W_{i,n+m}(\gamma_{n+m}),
\]
where
\[
\gamma_{n+m}(b_1,\ldots,b_{m+n})=
\alpha_n(\beta_m(b_1,\ldots,b_{m})b_{m+1},b_{m+2},\ldots,b_{m+n});
\]
\item
\begin{align*}
V_{i,n}(\alpha_n)W_{i,m}(\beta_m)&=V_{i,n+m-1}(\gamma_{n+m-1}), \\
W_{i,n}(\alpha_n)W_{i,m}(\beta_m)&=W_{i,n+m-1}(\gamma_{n+m-1}),
\end{align*}
where
\[
\gamma_{n+m-1}(b_1,\ldots,b_{m+n-1})=
\alpha_n(\beta_m(b_1,\ldots,b_{m}),b_{m+1},b_{m+2},\ldots,b_{m+n-1});
\]
\item
\[
\lambda(b)V_{i,n}(\alpha_n)=V_{i,n}(b\alpha_n),
\]
\item if $i'\ne i$ and $n\ge1$, then
\[
V_{i,n}(\alpha_n)L_{i'}=0=W_{i,n}(\alpha_n)L_{i'}.
\]
\end{enumerate}
\end{lemma}

\begin{prop}\label{prop:free}
For $i\in I$ let $\Afr_i\subseteq\Bc(\Fc)\cap\rho(B)'$ be the subalgebra
generated by
\[
\lambda(B)\cup\{L_i\}\cup\{V_{i,n}(\alpha_n)\mid n\in\Nats,\,\alpha_n\in\Bc_n(B)\}
\cup\{W_{i,n}(\alpha_n)\mid n\in\Nats,\,\alpha_n\in\Bc_n(B)\}.
\]
Then the family $(\Afr_i)_{i\in I}$ is free with respect to $\Ec$.
\end{prop}
\begin{proof}
Using Lemma~\ref{lem:rel},
we see that every element of $\Afr_i$ can be written as a sum of finitely
many terms of the following forms:
\renewcommand{\labelenumi}{(\roman{enumi})}
\begin{enumerate}
\item $\lambda(b)$
\item $\lambda(b_0)L_i\lambda(b_1)\cdots L_i\lambda(b_n)$
\item $V_{i,n}(\alpha_n)$
\item $\lambda(b_0)L_i\lambda(b_1)L_i\cdots \lambda(b_{k})L_iV_{i,n}(\alpha_n)$
\item $\lambda(b)W_{i,n}(\alpha_n)$
\item $\lambda(b_0)L_i\lambda(b_1)L_i\cdots \lambda(b_{k-1})L_i\lambda(b_k)W_{i,n}(\alpha_n)$.
\end{enumerate}
Now all terms of the forms (ii)--(vi) lie in $\ker\Ec$, while $\Ec(\lambda(b))=b$.
Therefore, $\Afr_i\cap\ker\Ec$ is the set of all finite sums of terms of the forms (ii)--(vi).

Let $p\in\Nats$ with $p\ge2$ and take $i_1,\ldots,i_p\in I$
with $i_1\ne i_2,\,i_2\ne i_3,\ldots,i_{p-1}\ne i_p$.
Suppose $a_j\in \Afr_{i_j}\cap\Ec$ ($1\le j\le p$) and let us show
$\Ec(a_1\cdots a_p)=0$.
From Lemma~\ref{lem:rel} part~(vi),
we see $a_1a_2\cdots a_p=0$ unless either $\forall j$ $a_j$ is of the form~(ii)
or $\forall j$ $a_j$ is of the form~(iii) or~(v).
But $V_{i,n}(\alpha_n)\Omega=0=W_{i,n}(\alpha_n)\Omega$ when $n\ge1$, so if $a_p$ is
of the form~(iii) or~(v), then $\Ec(a_1\cdots a_p)=0$.
We are left to consider the case when $a_1\ldots a_p$ can be written as
\begin{multline*}
(\lambda(b_0)L_{i_1}\lambda(b_1^{(1)})L_{i_1}\lambda(b_2^{(1)})\cdots L_{i_1}\lambda(b_{k(1)}^{(1)}))
(L_{i_2}\lambda(b_1^{(2)})\cdots L_{i_2}\lambda(b_{k(2)}^{(2)}))\cdots \\
\cdots(L_{i_p}\lambda(b_1^{(p)})\cdots L_{i_p}\lambda(b_{k(p)}^{(p)})),
\end{multline*}
where all $k(j)\ge1$.
But in this case, clearly $\Ec(a_1\cdots a_p)=0$.
\end{proof}

\begin{lemma}\label{lem:Nthmoment}
Let $N\in\Nats$ and for every $n\in\{0,1,\ldots,N\}$ let $\alpha_n\in\Bc_n(B)$.
Fix $i\in I$ and let
\begin{align*}
X&=\sum_{n=0}^{N-1}(V_{i,n}(\alpha_n)+W_{i,n}(\alpha_n)) \\[1ex]
Y&=X+V_{i,N}(\alpha_N)+W_{i,N}(\alpha_N).
\end{align*}
Then for any $b_0,\ldots,b_N\in B$, we have
\[
\Ec(b_0Yb_1Y\cdots b_NY)=b_0\alpha_N(b_1\alpha_0,b_2\alpha_0,\ldots,b_N\alpha_0)
+\Ec(b_0Xb_1X\cdots b_NX).
\]
\end{lemma}
\begin{proof}
To evaluate $\Ec(b_0Yb_1Y\cdots b_NY)$,
first write
\[
Y=\sum_{n=0}^N(V_{i,n}(\alpha_n)+W_{i,n}(\alpha_n))
\]
and distribute.
Now using the creation and annihilation properties
of the $W_{i,n}(\alpha_n)$ and $V_{i,n}(\alpha_n)$
operators, we see that the only term involving $\alpha_N$
to contribute a possibly nonzero quantity to $\Ec(b_0Yb_1Y\cdots b_NY)$
is
\[
\Ec(b_0V_{i,N}(\alpha_N)b_1W_{i,0}(\alpha_0)\cdots b_NW_{i,0}(\alpha_0)),
\]
whose value is $b_0\alpha_N(b_1\alpha_0,b_2\alpha_0,\ldots,b_N\alpha_0)$.
The other terms involve only $\alpha_0,\ldots,\alpha_{N-1}$ and
their sum is $\Ec(b_0Xb_1X\cdots b_NX)$.
\end{proof}

\begin{prop}\label{prop:Smom}
Let $(A,E)$ be a $B$--valued Banach noncommutative probability space and let $a\in A$, $N\in\Nats$.
Suppose $E(a)$ is an invertible element of $B$.
Let $\alpha_0=E(a)$.
Then there are $\alpha_1,\ldots,\alpha_N$, with $\alpha_n\in\Bc_n(B)$, such that if
\[
X=\sum_{n=0}^N(V_{i,n}(\alpha_n)+W_{i,n}(\alpha_n))\in\Bc(\Fc),
\]
then
\begin{equation}\label{eq:EcE}
\Ec(b_0Xb_1X\cdots b_kX)=E(b_0ab_1a\cdots b_ka)
\end{equation}
for all $k\in\{1,\ldots,N\}$ and all $b_0,\ldots,b_N\in B$.
\end{prop}
\begin{proof}
Using Lemma~\ref{lem:Nthmoment}, The maps $\alpha_k$ can be chosen
recursively in $k$ so that~\eqref{eq:EcE} holds.
\end{proof}

For the remainder of this section, we take $I=\{1,2\}$.
\begin{lemma}\label{lem:X}
Let $\alpha_0\in B$ be invertible.
Let $N\in\Nats$ and choose $\alpha_n\in\Bc_n(B)$ for $n\in\{1,\ldots,N\}$,
and let 
\[
F(b)=\alpha_0+\sum_{n=1}^N\alpha_n(b,\ldots,b).
\]
Note that $F(b)$ is invertible for $\|b\|$ sufficiently small.
Let
\begin{equation}\label{eq:X}
X=\sum_{n=0}^N(V_{1,n}(\alpha_n)+W_{1,n}(\alpha_n))\in\Bc(\Fc).
\end{equation}
Then the S--transform of $X$ is
$S_X(b)=F(b)^{-1}$.
\end{lemma}
\begin{proof}
For $b\in B$, $\|b\|<1$, let
\[
\omega_b=\Omega+\sum_{k=1}^\infty(b\delta_1)^{\otimes k}\otimes1\in\Fc
\]
We have $V_{1,0}(\alpha_0)\omega_b=\alpha_0\omega_b$ and, for $n\ge1$,
\[
V_{1,n}(\alpha_n)\omega_b=\alpha_n(b,\ldots,b)\Omega+\sum_{k=n+1}^\infty
\alpha_n(b,\ldots,b)(b\delta_1)^{\otimes(k-n)}\otimes1
=\alpha_n(b,\ldots,b)\omega_b.
\]
Moreover, $W_{1,0}(\alpha_0)\omega_b=\alpha_0L_1\omega_b$ and, for $n\ge1$,
\begin{align*}
W_{1,n}(\alpha_n)\omega_b&=\alpha_n(b,\ldots,b)\delta_1\otimes1
+\sum_{k=n+1}^\infty\alpha_n(b,\ldots,b)\delta_1\otimes(b\delta_1)^{\otimes(k-n)}\otimes1 \\
&=\alpha_n(b,\ldots,b)L_1\omega_b.
\end{align*}
Thus,
\[
X\omega_b=F(b)(1+L_1)\omega_b.
\]
For $\|b\|$ sufficiently small, we get
\begin{gather*}
F(b)^{-1}X\omega_b=\omega_b+L_1\omega_b \\
bF(b)^{-1}X\omega_b=b\omega_b+(\omega_b-\Omega) \\
\Omega=(1+b)\omega_b-bF(b)^{-1}X\omega_b \\
\Omega=(1-bF(b)^{-1}X(1+b)^{-1})(1+b)\omega_b \\
(1-bF(b)^{-1}X(1+b)^{-1})^{-1}\Omega=(1+b)\omega_b \\
\Ec((1-bF(b)^{-1}X(1+b)^{-1})^{-1})\begin{aligned}[t]
 =&P((1+b)\omega_b) \\
 =&1+b.\end{aligned}
\end{gather*}
Conjugating with $(1+b)$ yields
\[
1+b=\Ec((1-(1+b)^{-1}bF(b)^{-1}X)^{-1})=1+\Psi_X((1+b)^{-1}bF(b)^{-1}).
\]
Hence,
\[
\Psi_X\cinv(b)=(1+b)^{-1}bF(b)^{-1}
\]
and $S_X(b)=F(b)^{-1}$.
\end{proof}

\begin{lemma}\label{lem:XY}
Let $\alpha_0,\ldots,\alpha_n$, $F$ and $X$ be as in Lemma~\ref{lem:X}.
Let $\beta_0\in B$ be invertible and let $\beta_n\in\Bc_n(B)$ for $n\in\{1,\ldots,N\}$.
Let
\[
G(b)=\beta_0+\sum_{n=1}^N\beta_n(b,\ldots,b)
\]
and let
\begin{equation}\label{eq:Y}
Y=\sum_{n=0}^N(V_{2,n}(\alpha_n)+W_{2,n}(\alpha_n))\in\Bc(\Fc).
\end{equation}
Then the S--transform of $XY$ is
\begin{equation}\label{eq:SXY}
S_{XY}(b)=G(b)^{-1}F(G(b)bG(b)^{-1})^{-1}=S_Y(b)S_X(S_Y(b)^{-1}bS_Y(b)).
\end{equation}
\end{lemma}
\begin{proof}
From Lemma~\ref{lem:X}, we have $S_Y(b)=G(b)^{-1}$ and $S_X(b)=F(b)^{-1}$,
so the final equality in~\eqref{eq:SXY} is true.
For $b\in B$ let
\[
Z_b=bL_2+bG(b)^{-1}L_1G(b)+bG(b)^{-1}L_1G(b)L_2\in\Bc(\Fc)
\]
and insist that $\|b\|$ be so small that $\|Z_b\|<1$.
Let
\[
\sigma_b=(1-Z_b)^{-1}\Omega=\Omega+\sum_{k=1}^\infty Z_b^k\Omega.
\]
Using Lemma~\ref{lem:rel}, we find for $n,k\ge0$,
\[
V_{2,n}(\beta_n)Z_b^k=
\begin{cases}
V_{2,n-k}(\betat_{n-k}),&k<n, \\
\beta_n(b,\ldots,b),&k=n, \\
\beta_n(b,\ldots,b)Z_b^{k-n},&k>n
\end{cases}
\]
and
\[
W_{2,n}(\beta_n)Z_b^k=
\begin{cases}
W_{2,n-k}(\betat_{n-k}),&k<n, \\
\beta_n(b,\ldots,b)L_2,&k=n, \\
\beta_n(b,\ldots,b)L_2Z_b^{k-n},&k>n,
\end{cases}
\]
where
\[
\betat_{n-k}(b_1,\ldots,b_{n-k})=\beta_n(\underset{k}{\underbrace{b,\ldots,b}},b_1,\ldots,b_{n-k}).
\]
Therefore,
\[
V_{2,n}(\beta_n)Z_b^k\Omega=
\begin{cases}
0,&k<n, \\
\beta_n(b,\ldots,b)\Omega,&k=n, \\
\beta_n(b,\ldots,b)Z_b^{k-n}\Omega,&k>n
\end{cases}
\]
and
\[
W_{2,n}(\beta_n)Z_b^k\Omega=
\begin{cases}
0,&k<n, \\
\beta_n(b,\ldots,b)L_2\Omega,&k=n, \\
\beta_n(b,\ldots,b)L_2Z_b^{k-n}\Omega,&k>n
\end{cases}
\]
and we get
\[
Y\sigma_b=G(b)(1+L_2)\sigma_b.
\]
Letting $b'=G(b)bG(b)^{-1}$, we similarly find for $n,k\ge0$,
\[
V_{1,n}(\alpha_n)G(b)Z_b^k=
\begin{cases}
V_{1,n-k}(\alphat_{n-k})G(b),&k<n, \\
\alpha_n(b',\ldots,b')G(b)(1+L_2),&k=n, \\
\alpha_n(b',\ldots,b')G(b)(1+L_2)Z_b^{k-n},&k>n
\end{cases}
\]
and
\[
W_{1,n}(\alpha_n)G(b)Z_b^k=
\begin{cases}
W_{1,n-k}(\alphat_{n-k})G(b),&k<n, \\
\alpha_n(b',\ldots,b')L_1G(b)(1+L_2),&k=n, \\
\alpha_n(b',\ldots,b')L_1G(b)(1+L_2)Z_b^{k-n},&k>n,
\end{cases}
\]
where
\[
\alphat_{n-k}(b_1,\ldots,b_{n-k})=\alpha_n(\underset{k}{\underbrace{b',\ldots,b'}},b_1,\ldots,b_{n-k}).
\]
Therefore, we get
\[
XY\sigma_b=F(b')(1+L_1)G(b)(1+L_2)\sigma_b.
\]
Thus, for $\|b\|$ sufficiently small we get
\begin{gather*}
F(b')^{-1}XY=(1+L_1)G(b)(1+L_2)\sigma_b \\
F(b')^{-1}XY=G(b)\sigma_b+(G(b)L_2+L_1G(b)+L_1G(b)L_2)\sigma_b \\
bG(b)^{-1}F(b')^{-1}XY\sigma_b=b\sigma_b+Z_b\sigma_b \\
bG(b)^{-1}F(b')^{-1}XY\sigma_b=b\sigma_b+(\sigma_b-\Omega) \\
\Omega=((1+b)-bG(b)^{-1}F(b')^{-1}XY)\sigma_b \\
\Omega=(1-bG(b)^{-1}F(b')^{-1}XY(1+b)^{-1})(1+b)\sigma_b \\
(1-bG(b)^{-1}F(b')^{-1}XY(1+b)^{-1})^{-1}\Omega=(1+b)\sigma_b \\
\Ec((1-bG(b)^{-1}F(b')^{-1}XY(1+b)^{-1})^{-1})\begin{aligned}[t]
=&P((1+b)\sigma_b) \\
=&1+b.\end{aligned}
\end{gather*}
Conjugating with $(1+b)$ yields
\[
\Psi_{XY}((1+b)^{-1}bG(b)^{-1}F(b')^{-1})
=\Ec((1-(1+b)^{-1}bG(b)^{-1}F(b')^{-1}XY)^{-1})-1=b.
\]
Hence,
\[
\Psi_{XY}\cinv(b)=(1+b)^{-1}bG(b)^{-1}F(b')^{-1}
\]
and~\eqref{eq:SXY} holds.
\end{proof}

\begin{proof}[Proof of Theorem~\ref{thm:main}]
The formula~\eqref{eq:main} asserts the equality of the germs
of two analytic $B$--valued functions.
This is equivalent to asserting the equality of the $n$th terms
in their respective power series expansions around zero, for every $n\ge0$.
By Lemmas~\ref{lem:FG} and~\ref{lem:Finv}, the $n$th term, call it RHS$_n$, in
the expansion for the right hand side of~\eqref{eq:main} depends only
on the $0$th through the $n$th terms of the power series expansions
for $S_x(b)$ and $S_y(b)$.
Hence, by Proposition~\ref{prop:nth}, RHS$_n$ depends only on the
moment functions $\mu_{x,1},\,\ldots,\,\mu_{x,n}$ and
$\mu_{y,1},\,\ldots,\,\mu_{y,n}$.
On the other hand, again by Proposition~\ref{prop:nth}, the $n$th term in the
power series expansion for the left hand side of~\eqref{eq:main},
call it LHS$_n$, depends only on
$\mu_{xy,1},\,\ldots,\,\mu_{xy,n}$.
But by freeness of $x$ and $y$, for each $k\ge1$ the moment function
$\mu_{xy,k}$ depends only on $\mu_{x,1},\,\ldots,\,\mu_{x,k}$
and $\mu_{y,1},\,\ldots,\,\mu_{y,k}$.
Thus, both LHS$_n$ and RHS$_n$ depend only on
$\mu_{x,1},\,\ldots,\,\mu_{x,n}$ and
$\mu_{y,1},\,\ldots,\,\mu_{y,n}$.

Hence, in order to prove~\eqref{eq:main} at the level of the $n$th terms
in the power series expansion, it will suffice to prove~\eqref{eq:main}
for some free pair $X$ and $Y$ of elements in a Banach noncommutative probability space over $B$,
whose first $n$ moment functions agree
with those of $x$ and $y$, respectively.
However, by Propositions~\ref{prop:free} and~\ref{prop:Smom},
such $X$ and $Y$ can be chosen of the forms~\eqref{eq:X} and~\eqref{eq:Y}.
By Lemma~\ref{lem:XY}, the equality~\eqref{eq:main} holds for these operators.
\end{proof}

\section{A proof of the additivity of the R--transform over a Banach algebra}
\label{sec:R-t}

The R--transform over a general unital algebra $B$ has been well understood
since Voiculescu's work~\cite{V95} (and see also Speicher's approach in~\cite{Sp}). 
However, for completeness,
in this section we offer a new proof, using the techniques
and constructions of the previous two sections, of the additivity of the R--transform
for free random variables in a Banach noncommutative probability space.
This proof is, of course, analogous to Haagerup's proof of Theorem 2.2 of~\cite{H}
in the scalar--valued case.

Let $(A,E)$ be a Banach noncommutative probability space over $B$
and let $a\in A$.
Consider the function
\[
\Cc_a(b)=E((1-ba)^{-1}b)=\sum_{n=0}^\infty E((ba)^nb),
\]
defined and analytic for $\|b\|<\|a\|^{-1}$.
We have $\Cc_a(b)=b+b\Phi_a(b)b$, where $\Phi_a$ is as in~\eqref{eq:Phi}.
Since the Fr\'echet differential of $\Cc_a$ at $b=0$ is the identity map, $\Cc_a$
is invertible with respect to composition in a neighborhood of zero.
\begin{prop}\label{prop:R}
There is a unique $B$--valued analytic function $R_a$, defined in a neighborhood of $0$ in $B$,
such that
\begin{equation}\label{eq:CR}
\Cc_a\cinv(b)=(1+bR_a(b))^{-1}b=b(1+R_a(b)b)^{-1}.
\end{equation}
\end{prop}
\begin{proof}
Again, uniqueness is clear by the power series expansions.

The right--most equality in~\eqref{eq:CR} holds for any analytic function $R_a$.
We seek a function $R_a$ such that
\[
\Cc_a((1+bR_a(b))^{-1}b)=b.
\]
But
\[
\Cc_a((1+bR_a(b))^{-1}b)=\begin{aligned}[t]
  &(1+bR_a(b))^{-1}b \\
  &+(1+bR_a(b))^{-1}b\,\Phi_a\big((1+bR_a(b))^{-1}b\big)\,(1+bR_a(b))^{-1}b,\end{aligned}
\]
so it will suffice to find $R_a$ so that any of the following hold:
\begin{gather}
\notag
(1+bR_a(b))^{-1}+(1+bR_a(b))^{-1}b\,\Phi_a\big((1+bR_a(b))^{-1}b\big)\,(1+bR_a(b))^{-1}=1, \\
\notag
1+b\,\Phi_a\big((1+bR_a(b))^{-1}b\big)\,(1+bR_a(b))^{-1}=1+bR_a(b), \\
\notag
b\,\Phi_a\big((1+bR_a(b))^{-1}b\big)\,(1+bR_a(b))^{-1}=bR_a(b), \\
\Phi_a\big((1+bR_a(b))^{-1}b\big)\,(1+bR_a(b))^{-1}=R_a(b). \label{eq:Ra}
\end{gather}
However, $R_a(0)=E(a)$ is a solution of~\eqref{eq:Ra} at $b=0$, and the
Fr\'echet differential of the function $x\mapsto \Phi_a((1+bx)^{-1}b)(1+bx)^{-1}-x$ at $b=0$
is the negative of the identity map, hence is invertible.
The implicit function theorem of Hildebrandt and Graves~\cite{HG}
(see also the discussion on p.\ 655 of~\cite{G}) guarantees the existence of $R_a$.
\end{proof}

The {\em R--transform} of $a$ is defined to be the analytic function $R_a$ from
Proposition~\ref{prop:R}.

Analogously to Proposition~\ref{prop:nth}, we have the following.
\begin{prop}\label{prop:Rn+1}
The $n$th term $(R_a)_n(b,\ldots,b)$ in the power series expansion for $R_a$ about zero 
depends only on the first $n+1$ moment functions $\mu_{a,1},\ldots,\mu_{a,n+1}$ of $a$.
\end{prop}

Here is the analogue to Lemma~\ref{lem:Nthmoment}, which can be proved similarly.
\begin{lemma}
Let $N\in\Nats$ and for every $n\in\{0,1,\ldots,N\}$ let $\alpha_n\in\Bc_n(B)$.
Fix $i\in I$ and let
\begin{align*}
X&=L_i+\sum_{n=0}^{N-1}V_{i,n}(\alpha_n) \\[1ex]
Y&=X+V_{i,N}(\alpha_N).
\end{align*}
Then for any $b_0,\ldots,b_N\in B$, we have
\[
\Ec(b_0Yb_1Y\cdots b_NY)=b_0\alpha_N(b_1,b_2,\ldots,b_N)
+\Ec(b_0Xb_1X\cdots b_NX).
\]
\end{lemma}

We immediately get the following analogue of Proposition~\ref{prop:Smom}.
\begin{prop}\label{prop:Rmom}
Let $(A,E)$ be a $B$--valued Banach noncommutative probability space and let $a\in A$, $N\in\Nats$.
Then there are $\alpha_0,\alpha_1,\ldots,\alpha_N$, with $\alpha_n\in\Bc_n(B)$, such that if
\[
X=L_i+\sum_{n=0}^NV_{i,n}(\alpha_n)\in\Bc(\Fc),
\]
then
\[
\Ec(b_0Xb_1X\cdots b_kX)=E(b_0ab_1a\cdots b_ka)
\]
for all $k\in\{1,\ldots,N\}$ and all $b_0,\ldots,b_N\in B$.
\end{prop}

Now we have the following analogues of Lemmas~\ref{lem:X} and~\ref{lem:XY}.

\begin{lemma}\label{lem:R-tX}
Let $N\in\Nats$ and choose $\alpha_n\in\Bc_n(B)$ for $n\in\{0,1,\ldots,N\}$,
and let 
\[
F(b)=\alpha_0+\sum_{n=1}^N\alpha_n(b,\ldots,b).
\]
Let
\[
X=L_1+\sum_{n=0}^NV_{1,n}(\alpha_n)\in\Bc(\Fc).
\]
Then the R--transform of $X$ is
$R_X(b)=F(b)$.
\end{lemma}
\begin{proof}
With $\omega_b$ defined as in the proof of Lemma~\ref{lem:X},
we have
\begin{gather*}
X\omega_b=L_1\omega_b+F(b)\omega_b \\
bX\omega_b=(\omega_b-\Omega)+bF(b)\omega_b \\
(1+bF(b)-bX)\omega_b=\Omega \\
(1-bX(1+bF(b))^{-1})^{-1}\Omega=(1+bF(b))\omega_b \\
\Ec((1-bX(1+bF(b))^{-1})^{-1})\begin{aligned}[t]&=P((1+bF(b))\omega_b) \\
  &=1+bF(b).\end{aligned}
\end{gather*}
Conjugating yields 
\[
\Ec((1-(1+bF(b))^{-1}bX)^{-1})=1+bF(b),
\]
so
\[
\Cc_X((1+bF(b))^{-1}b)=\Ec((1-(1+bF(b))^{-1}bX)^{-1})(1+bF(b))^{-1}b=b.
\]
Thus,
\[
\Cc_X\cinv(b)=(1+bF(b))^{-1}b
\]
and $R_a(b)=F(b)$.
\end{proof}

\begin{lemma}\label{lem:R-tXY}
Let $\alpha_0,\ldots,\alpha_n$, $F$ and $X$ be as in Lemma~\ref{lem:R-tX}.
Let $\beta_n\in\Bc_n(B)$ for $n\in\{0,1,\ldots,N\}$.
Let
\[
G(b)=\beta_0+\sum_{n=1}^N\beta_n(b,\ldots,b)
\]
and let
\[
Y=L_2+\sum_{n=0}^NV_{2,n}(\alpha_n)\in\Bc(\Fc).
\]
Then the R--transform of $X+Y$ is
\[
R_{X+Y}(b)=F(b)+G(b)=R_X(b)+R_Y(b).
\]
\end{lemma}
\begin{proof}
For $b\in B$ with $\|b\|<1/2$, let
\[
\sigma_b=(1-b(L_1+L_2))^{-1}\Omega
=\Omega+\sum_{k=1}^\infty(b\delta_1+b\delta_2)^{\otimes k}\otimes 1\in\Fc.
\]
Then
\begin{gather*}
(X+Y)\sigma_b=(L_1+L_2)\sigma_b+(F(b)+G(b))\sigma_b \\
b(X+Y)\sigma_b=(\sigma_b-\Omega)+b\,(F(b)+G(b))\sigma_b.
\end{gather*}
Now arguing as in the proof of Lemma~\ref{lem:R-tX} above yields $R_{X+Y}(b)=F(b)+G(b)$.
\end{proof}

Finally, we get a proof, which is analogous to our proof of Theorem~\ref{thm:main}, of the additivity
of the R--transform in a Banach noncommutative probability space.
\begin{thm}[\cite{V95}]
Let $B$ be a unital complex Banach algebra and let $(A,E)$ be a $B$--valued
Banach noncommutative probability space.
Let $x,y\in A$ be free in $(A,E)$.
Then
\[
R_{x+y}(b)=R_x(b)+R_y(b).
\]
\end{thm}
\begin{proof}
By Proposition~\ref{prop:Rn+1}, it will suffice to show
that given $n\in\Nats$ we have $R_{X+Y}=R_X+R_Y$ for
some free pair $X$ and $Y$ of elements in a Banach noncommutative probability space over $B$
whose first $n$ moment functions agree with those of $x$ and $y$, respetively.
Precisely this fact follows from Proposition~\ref{prop:Rmom}, Proposition~\ref{prop:free}
and Lemmas~\ref{lem:R-tX} and~\ref{lem:R-tXY}.
\end{proof}

\bibliographystyle{plain}

\begin{thebibliography}{10}

\bibitem{Aa} L.\ Aagaard,
`A Banach algebra approach to amalgamated R-- and S--transforms,'
preprint (2004).

\bibitem{G} L.M. Graves,
`Topics in functional calculus,'
{\em Bull. Amer. Math. Soc.} {\bf 41} (1935), 641-662.
Correction, {\em ibid.} {\bf 42} (1936), 381-382.

\bibitem{H} U.\ Haagerup,
`On Voiculescu's R-- and S-- transforms for free non--commuting random variables,'
{\em Free Probability Theory,} D.~Voiculescu, (Ed.),
Fields Inst.\ Commun.\ {\bf 12} (1997), 127-148.

\bibitem{HG} T.H.\ Hildebrandt and L.M.\ Graves,
`Implicit functions and their differentials in general analysis,'
{\em Trans.\ Amer.\ Math.\ Soc.\ } {\bf 29} (1927), 127-153.

\bibitem{HP} E.\ Hille and R.S.\ Phillips,
{\em Functional Analysis and Semi--groups,}
revised edition, American Mathematical Society, Providence,
1957.

\bibitem{Sp} R.\ Speicher,
`Combinatorial theory of the free product with amalgamation and operator-valued free probability theory,'
{\em Mem. Amer. Math. Soc.}, {\bf 132} 1998, no.\ 627.

\bibitem{V85} D.\ Voiculescu,
`Symmetries of some reduced free product $C^*$-algebras,'
{\em Operator Algebras and Their Connections with Topology and Ergodic Theory},
H. Araki, C.C.\ Moore, \c S.\ Str\u atil\u a and D.\ Voiculescu, (Eds.),
Lecture Notes in Mathematics, Volume 1132, Springer-Verlag, 1985, 556--588.

\bibitem{V87} D.\ Voiculescu,
`Multiplication of certain noncommuting random variables,'
{\em J. Operator Theory} {\bf 18} (1987), 223-235.

\bibitem{V95} D.\ Voiculescu,
`Operations on certain non-commutative operator--valued random variables,'
{\em Recent Advances in Operator Algebras (Orl\'eans, 1992)},
Ast\'erisque  No. 232 (1995), pp.\ 243-275.

\bibitem{VDN} D.V.\ Voiculescu, K.J.~Dykema, A.~Nica,
{\em Free Random Variables},
CRM Monograph Series {\bf 1}, American Mathematical Society, 1992.

\end{thebibliography}

\end{document}